\documentclass{elsart3-1}

\usepackage{amsmath,amssymb,graphicx,epsfig,color}
\usepackage[english,francais]{babel}

\newtheorem{e-proposition}[theorem]{Proposition}

\newtheorem{e-definition}[theorem]{Definition\rm}

\renewcommand{\theequation}{\arabic{equation}}
\def\zz{\mathbb{Z}}
\def\rr{\mathbb{R}}

\def\og{\leavevmode\raise.3ex\hbox{$\scriptscriptstyle\langle\!\langle$~}}
\def\fg{\leavevmode\raise.3ex\hbox{~$\!\scriptscriptstyle\,\rangle\!\rangle$}}

\begin{document}

\begin{frontmatter}

\selectlanguage{english}
\title{High frequency wave packets for the Schr\"{o}dinger equation and its numerical approximations}

\vspace{-2.6cm}

\selectlanguage{francais}
\title{Paquets d'ondes \`{a} haute fr\'{e}quence pour l'\'{e}quation de Schr\"{o}dinger et ses approximations num\'{e}riques}

\selectlanguage{english}
\author[BCAM]{Aurora Marica},
\ead{marica@bcamath.org}
\author[Iker,BCAM]{Enrique Zuazua}
\ead{zuazua@bcamath.org}

\address[Iker]{Ikerbasque, Basque Foundation for Science, Alameda Urquijo 36-5, Plaza Bizkaia, 48011, Bilbao, Basque Country, Spain}

\address[BCAM]{BCAM - Basque Center for Applied Mathematics, Bizkaia Technology Park 500, 48160, Derio, Basque Country, Spain}

\begin{abstract}
We build Gaussian wave packets for the linear Schr\"{o}dinger equation and its finite difference space semi-discretization and illustrate the lack of uniform dispersive properties of the numerical solutions as established in \cite{LivEZDispSchr}. It is by now well known that  bigrid algorithms provide filtering mechanisms allowing to recover the uniformity of the dispersive properties as the mesh size goes to zero. We analyze and illustrate numerically how these high frequency wave packets split and propagate under these bigrid filtering mechanisms, depending on how the fine grid/coarse grid filtering is implemented.
\vskip 0.5\baselineskip

\selectlanguage{francais}
\noindent{\bf R\'esum\'e
 } \vskip 0.5\baselineskip \noindent
On construit des paquets d'ondes gaussiennes pour l'\'{e}quation de Schr\"{o}dinger lin\'{e}aire continue unidimensionnelle ainsi que pour sa semi-discr\'{e}tisation en espace par diff\'{e}rences finies. On illustre num\'{e}riquement le manque d'uniformit\'{e} par rapport au pas du maillage des propri\'{e}t\'{e}s de dispersion des solutions num\'{e}riques d\'{e}montr\'{e} dans \cite{LivEZDispSchr}. Par ailleurs, il est bien connu que les algorithmes bi-grilles sont des m\'{e}canismes de filtrage efficaces pour r\'{e}cup\'{e}rer l'uniformit\'{e} des propri\'{e}t\'{e}s dispersives. On analyse la fa\c{c}on dont les solutions bi-grilles correspondant \`{a} plusieurs projections de la grille fine sur la grossi\`{e}re, se divisent en plusieurs paquets d'ondes, chacun se propageant diff\'{e}rement. On repr\'{e}sente num\'{e}riquement ces ph\'{e}nom\`{e}nes et on montre que ce comportement est en accord avec les r\'{e}sultats th\'{e}oriques connus sur la dispersion des solutions bi-grilles.
\end{abstract}
\end{frontmatter}

\selectlanguage{francais}
\section*{Version fran\c{c}aise abr\'eg\'ee}
Dans \cite{LivEZDispSchr}, il a \'{e}t\'{e} d\'{e}montr\'{e} que les solutions de l'\'{e}quation de Schr\"{o}dinger lin\'{e}aire unidimensionnelle semi-discretis\'{e}e en espace par dif\-f\'{e}\-rences finies sur une grille uniforme ne v\'{e}rifient pas les estimations dispersives de l'\'{e}quation de Schr\"{o}dinger continue (estimations $L^q_t-L^p_x$ et gain d'une demi-d\'{e}riv\'{e}e localement en espace) uniform\'{e}ment par rapport au pas du maillage. Ce comportement pathologique est d\^{u} \`{a} l'existence de points critiques qui annulent l'une des deux premi\`{e}res d\'{e}riv\'{e}es de la relation de dispersion.

Une mani\`{e}re de r\'{e}cup\'{e}rer l'uniformit\'{e} de ces estimations consiste \`{a} faire un filtrage bi-grille des donn\'{e}es initiales. Ce m\'{e}canisme a \'{e}t\'{e} introduit dans le travail pionnier de R. Glowinski \cite{GloStokesWaves} sur le filtrage bi-grille pour le contr\^{o}le des ondes. Son efficacit\'{e} a \'{e}t\'{e} d\'{e}montr\'{e}e rigureusement dans le contexte des propri\'{e}t\'{e}s dispersives  de l'\'{e}quation de Schr\"{o}dinger dans \cite{LivEZDispSchr}. Plus pr\'{e}cisement, tout en gardant le sch\'{e}ma d'approximation semi-discret conservatif en diff\'{e}rences finis, il s'agit de consid\'{e}rer uniquement des donn\'{e}es i\-ni\-ti\-ales dans un maillage quatre fois plus grossier et de les prolonger par interpolation lin\'{e}aire sur la grille fine.

Notre objectif principal dans cette Note est d'illustrer num\'{e}riquement les effets pathologiques des solutions num\'{e}riques et d'\'{e}tudier dans quelle mesure ils disparaissent lorsqu'on applique les m\'{e}canismes de filtrage bi-grille.

Nous consid\'{e}rons des donn\'{e}es initiales oscillatoires obtenues \`{a} partir d'un profil gaussien concentr\'{e} en Fourier autour d'un certain nombre d'onde. Afin de comparer les diff\'{e}rentes solutions num\'{e}riques obtenues par le filtrage bi-grille, avant d'utiliser l'interpolation lin\'{e}aire pour passer de la maille grossi\`{e}re \`{a} la fine, on doit projeter les donn\'{e}es initiales du maillage fin sur le grossier. On le fait de deux mani\`{e}res: d'une part, par restriction des fonctions discr\`{e}tes aux points du maillage grossier et, d'autre part, en moyennant les valeurs correspondant aux points de la grille fine situ\'{e}s dans le voisinage de chaque point dans la grille grossi\`{e}re. La premi\`{e}re peut g\'{e}n\'{e}rer, entre autres, des solutions qui n'oscillent pas et ne se propagent pas, concentr\'{e}es en Fourier autour du nombre d'onde $\xi=0$, donc, r\'{e}guli\`{e}res. La seconde peut donner lieu \`{a} des solutions tendant vers z\'{e}ro quand le pas du maillage tend vers z\'{e}ro.

Dans les deux cas, on constate que ce comportement est compatible avec les propri\'{e}t\'{e}s dispersives prouv\'{e}es dans \cite{LivEZDispSchr} pour la m\'{e}thode de diff\'{e}rences finies apr\`{e}s le filtrage bimaille.

\selectlanguage{english}
\textbf{Problem formulation.} Let us consider the $1-d$ linear continuous Schr\"{o}dinger equation (CSE):
\begin{equation}i\partial_tu(x,t)+\partial_x^2u(x,t)=0,\ x\in\rr,t\in\rr\setminus\{0\},\quad u(x,0)=\varphi(x),x\in\rr.\label{ContSchrEqn}\end{equation}
The solution of (\ref{ContSchrEqn}) is given by $u(x,t)=S(t)\varphi(x),$ where $S(t)=\exp(it\partial_x^2)$ is the Schr\"{o}dinger semigroup defined as $S(t)\varphi(x)=(G(\cdot,t)*\varphi)(x) $ and $G(x,t)$ is the fundamental solution of (\ref{ContSchrEqn}).

The solution of (\ref{ContSchrEqn}) verifies two \textit{dispersive properties}: the \textit{gain on the integrability} and the \textit{local smoothing effect}
\begin{equation}\|u\|_{L^q_t(\rr,L^p_x(\rr))}\leq c(p)\|\varphi\|_{L^2(\rr)}, \mbox{ respectively } \sup_R\Big(\frac{1}{R}\int_{\rr}\int_{-R}^R|\partial_x^{1/2}u(x,t)|^2\,dx\,dt\Big)^{1/2}\leq c\|\varphi\|_{L^2(\rr)},\label{GlobalStrichartz}\end{equation}
for $(p,q)$ such that the \textit{admissibility conditions} $2\leq p\leq\infty$ and $2/q=1/2-1/p$ hold. For detailed proofs and higher dimensional versions of these dispersive estimates, see \cite{LinaresPonce}. They play a key role in the proof of the well-posedness of the non-linear Schr\"{o}dinger equation in $L^2(\rr)$ or $H^1(\rr)$ (\cite{CazenaveSemilin}, \cite{LivEZDispSchr}).

On an uniform grid of size $h>0$ of the real line, $\mathcal{G}_h=\{x_j=jh: j\in\zz\}$, we introduce the discrete Laplacian of the sequence $\overrightarrow{f}^h=(f_j)_{j\in\zz}$ to be $\partial_h^2f_j=h^{-2}(f_{j+1}-2f_j+f_{j-1})$, and analyze the finite difference semi-discrete Schr\"{o}dinger equation (DSE):
\begin{equation}i\partial_tu_j(t)+\partial_h^2u_j(t)=0,\ j\in\zz,t\in\rr\setminus\{0\},\quad u_j(0)=\varphi_j,\ j\in\zz.\label{DiscSchrEqn}\end{equation}

Set $\Pi_h:=[-\pi/h,\pi/h]$. The solution of (\ref{DiscSchrEqn}) can be expressed by means of the inverse semi-discrete Fourier transform (SDFT) (cf. \cite{LivEZDispSchr}) as $u_j(t)=[\exp(it\partial_h^2)\overrightarrow{\varphi}^h]_j=\frac{1}{2\pi}\int_{\Pi_h}\widehat{\varphi}^h(\xi)\exp(itp_h(\xi))\exp(i\xi x_j)\,d\xi$. The symbol $p_h:\Pi_h\to\rr$ is defined as $p_h(\xi)=4h^{-2}\sin^2(\xi h/2)$ and $\widehat{\varphi}^h$ is the SDFT at scale $h$ of the initial data $\overrightarrow{\varphi}^h=(\varphi_j)_{j\in\zz}$.

In the continuous case, the Fourier symbol of the Laplacian is $p(\xi)=|\xi|^2$.  Its first order derivative, the so-called \textit{group velocity}, $\partial_{\xi}p(\xi)=2\xi$, vanishes only at $\xi=0$ and its second-order derivative, the so-called \textit{group acceleration}, $\partial_{\xi}^2p(\xi)=2$, does not vanish at any wave number. For the semi-discrete case, the following two pathologies of the symbol $p_h(\xi)$ were observed (cf. \cite{LivEZDispSchr}):

\begin{itemize}\item[p1.] The group velocity, $\partial_{\xi}p_h(\xi)=2h^{-1}\sin(\xi h)$, vanishes at $\xi=0$, but also at $\xi=\pm\pi/h$.

\item[p2.] The group acceleration, $\partial_{\xi}^2p_h(\xi)=2\cos(\xi h)$, vanishes at $\xi=\pm\pi/2h$.\end{itemize}

The pathologies (p1) and (p2) lead to the lack of uniform dispersive properties as $h\to 0$ and, more precisely, to the lack of local smoothing effects and of discrete Strichartz estimates, respectively (cf. \cite{LivEZDispSchr}).

Our main goal is to illustrate the effects of these pathologies on the numerical solutions and to analyze to which extent they disappear when applying the two-grid filtering mechanisms proposed in \cite{LivEZDispSchr} and inspired by the pioneering work by R. Glowinski \cite{GloStokesWaves} on the bi-grid filtering for the control of waves.

As proved in \cite{LivEZDispSchr}, an efficient mechanism to recover the uniformity as $h\to0$ of the dispersive properties is the \textit{bi-grid algorithm} introduced in \cite{GloStokesWaves}. This consists in solving the DSE on the fine grid of size $h$ with slow initial data obtained by linear interpolation from data given on a coarser grid of size $nh$. We will take $n=2^k$, with $k\geq 1$. For ratios $1/2^k$, $k\geq 2$, between the two grids (the fine one of size $h$ and the coarse one of size $2^kh$), both pathologies (p1) and (p2) are canceled and the uniform dispersivity is recovered (cf. \cite{LivEZDispSchr}).

Define the extension operator $\Gamma_k:\ell^2(2^kh\zz)\to\ell^2(h\zz)$ from the grid of size $2^kh$ to the one of size $h$, by linear interpolation, as follows: $(\Gamma_k f)_{2^kj+r}=(2^k-r)/2^k f_{2^kj}+r/2^k f_{2^kj+2^k},$ for all $j\in\zz$ and $0\leq r\leq 2^k-1$. If $\overrightarrow{f}\in\ell^2(2^kh\zz)$, the Fourier transform of this extension can be written as $\widehat{\Gamma_k f}^{h}(\xi)=b_k(\xi h)\widehat{f}^{2^kh}(\xi)$, for all $\xi\in\Pi_h,$ with weights $b_k(\eta)=\prod_{j=1}^k\cos^2(2^{j-2}\eta)$ vanishing quadratically at $\eta=\pm j\pi/2^{k-1}$, for all $1\leq j\leq 2^{k-1}$. The SDFT $\widehat{f}^{2^kh}$ is defined for $\xi\in\Pi_{2^kh}$, but for simplicity, we also denote in this manner its extension to $\Pi_h$ by $\pi/(2^{k-1}h)$-periodicity.

According to \cite{LivEZDispSchr}, when $k=2$, i. e. when the ratio between the two meshes is $1/4$, and the numerical scheme is restricted to this class of filtered initial data, the dispersive properties turn out to be uniform as $h \to 0$. More precisely, there exist two constants $C(p),C>0$ independent of $h$ such that the solution $\overrightarrow{u}^h(t)=\exp(it\partial_h^2)\Gamma_2\overrightarrow{\varphi}^h$ verifies
$$\|\overrightarrow{u}^h(t)\|_{L^q(\rr,\ell^p(h\zz))}\leq C(p)\|\Gamma_2\overrightarrow{\varphi}^h\|_{\ell^2(h\zz)}\mbox{ and }\sup_{R>0}\frac{1}{R}\int_{\rr}h\sum_{|x_j|\leq R}|\partial_h^{1/2}u_j(t)|^2\,dt\leq C\|\Gamma_2\overrightarrow{\varphi}^h\|_{\ell^2(h\zz)}^2.$$
Here, $\partial_h^s$ denotes the discrete fractional derivative, i.e.  $\partial_h^sf_j=(2\pi)^{-1}\int_{\Pi_h}p_h^{s/2}(\xi)\widehat{f}^h(\xi)\exp(i\xi x_j)\,d\xi$, with $s\geq 0$, and the discrete $\ell^p(h\zz)$-spaces are defined as usual (see, e.g. \cite{LivEZDispSchr}).

But, in practice, the initial data for the DSE are given on the fine grid of size $h$ as an approximation of the initial datum of the continuous Schr\"{o}dinger equation on the nodal points $x_j=jh$. Thus, they need, first, to be projected into the coarse one. We analyze two different projection operators from $\mathcal{G}_h$ to $\mathcal{G}_{2^kh}$:
\begin{equation}(\Lambda_k^rf)_{2^kj}=f_{2^kj}\mbox{ and } (\Lambda_k^af)_{2^kj}=\sum_{r=0}^{2^k-1}\Big(\frac{2^k-r}{2^{2k}}f_{2^kj+r}+\frac{r}{2^{2k}}f_{2^kj+r-2^k}\Big).\label{projectionsbigrid}\end{equation}
The superscripts $r$ and $a$ stand for \textit{restriction} and \textit{average}, the two key mechanisms on which these projections are based. More precisely, the projection $\Lambda_k^r$ restricts the function $(f_j)_{j\in\zz}$ on the fine grid to those $j$-s that are integer multiples of $2^k$. The projection $\Lambda_k^{a}$ takes as value at the point $x_{2^kj}$ an average of the values at the $2^{k+1}-1$ points surrounding $x_{2^kj}$ in the fine grid.

These projection operators can be represented in the Fourier space, $\Pi_{2^kh}$, as follows:
\begin{equation}\widehat{\Lambda_k^rf}^{2^kh}=\sum\limits_{j=-2^{k-1}}^{2^{k-1}-1}\widehat{f}^h\Big(\cdot+\frac{2j\pi}{2^kh}\Big)\mbox{ and }
\widehat{\Lambda_k^af}^{2^kh}=\sum\limits_{j=-2^{k-1}}^{2^{k-1}-1}\widehat{f}^h\Big(\cdot+\frac{2j\pi}{2^k h}\Big)b_k\Big(\cdot h+\frac{2j\pi}{2^k}\Big).
\label{projectionsbigridFour}\end{equation}

\smallskip\noindent\textbf{Behavior of solutions under quadratic dispersion relations and Gaussian initial data.} Our aim is to describe the behavior of the numerical solutions that the bi-grid algorithm produces when the initial data on the fine grid is a highly concentrated Gaussian wave packet.  For this, we introduce the Gaussian profile $\widehat{\sigma}_{\gamma}(\xi)=\sqrt{2\pi/\gamma}\exp(-|\xi|^2/(2\gamma)),$ where
\begin{equation}\gamma=\gamma(h)\mbox{ such that }\gamma h^{2/3}<<1\mbox{ and }\gamma>>1.\label{scalegamma}\end{equation}

Consider $\eta_0,\eta_1<\eta_2\in[-\pi,\pi]$ independent of $h$ and the following function depending on $\eta_0,\eta_1,\eta_2,\gamma$:
\begin{equation}\widehat{\varphi}_{\eta_0,\gamma}^{\eta_1,\eta_2}(\xi)=\widehat{\sigma}_{\gamma}(\xi-\eta_0/h)
\chi_{[\eta_1/h,\eta_2/h]}(\xi)\label{initialdataeta0}\end{equation}
For $\eta_0\in(-\pi,\pi]$, $\widehat{\varphi}_{\eta_0}$ represents a Gaussian profile concentrated around $\eta_0/h$ truncated to $[\eta_1/h,\eta_2/h]\subseteq\Pi_h$. Note that even for the CSE, for an easier comparison with the DSE,  we consider the same initial data supported in $\Pi_h$ in the Fourier space.

Taking into account that the dispersion relation of the CSE is quadratic, one can obtain an explicit representation formula for the corresponding solution $u_{\eta_0,\gamma}^{\eta_1,\eta_2}$. The same can be done, for more general dispersion relations, as it is for instance the case for the DSE, by taking the second order Taylor expansion. In the purely quadratic case we have:
\begin{equation}q_{\eta_0,h}(\xi)=\frac{1}{h^2}q_{\eta_0}(\xi h), \mbox{ with }q_{\eta_0}(\eta)\sim q_{\eta_0}(\eta_0)+q_{\eta_0}'(\eta_0)(\eta-\eta_0)+\frac{1}{2}q_{\eta_0}''(\eta_0)(\eta-\eta_0)^2.\label{quadraticdispersiongeneral}\end{equation}
This, together with the expression of the initial data (\ref{initialdataeta0}), allows to conclude that:
\begin{itemize}\item[i.] \textit{Propagation:} $u_{\eta_0,\gamma}^{\eta_1,\eta_2}$ propagates along the curve $x(t)=x^*-tq_{\eta_0}'(\eta_0)/h$.
\item[ii.] \textit{Time evolution of the support:} it expands like $(\gamma^{-1}+t^2\gamma (q_{\eta_0}''(\eta_0))^2)^{1/2}$.
\item[iii.] \textit{Time evolution of the amplitude:} it behaves like $c(1+t^2\gamma^2(q_{\eta_0}''(\eta_0))^2)^{-1/4}$, where $c=1$ if $\eta_0\not=\eta_1,\eta_2$ and $c=1/2$ if $\eta_0=\eta_1$ or $\eta_0=\eta_2$. For $q_{\eta_0}''(\eta_0)\not=0$, this large time behavior agrees with the $t^{-1/2}$-decay rate of the solutions of the CSE with $L^1$-data (see \cite{LivEZDispSchr}).
\end{itemize}

We now consider the CSE with two particular choices of initial data $\varphi$ given in the Fourier space as follows:
\begin{equation}\widehat{\varphi}_{\pi}(\xi)=\widehat{\varphi}^{-\pi,0}_{-\pi,\gamma}(\xi)+\widehat{\varphi}^{0,\pi}_{\pi,\gamma}(\xi)\mbox{ or }\widehat{\varphi}_{\eta_0}(\xi)=\widehat{\varphi}^{-\pi,\pi}_{\eta_0,\gamma}(\xi),\eta_0\in(-\pi,\pi).\label{initialdatawithoutbigrid}\end{equation}
In (\ref{initialdatawithoutbigrid}), we roughly consider initial data as in (\ref{initialdataeta0}), but with a particular choice of  $\eta_0$, $\eta_1$ and $\eta_2$ and with two superposed Gaussian profiles when they are supported in any of the end points $\pm\pi/h$. The initial data $\widehat{\varphi}_{\pi}$, having two picks, seems to concentrate around two wave numbers, $\pm\pi/h$, and this occurs in the continuous setting. But in the discrete one, since the two picks centered at $\pm\pi/h$ are located exactly on the boundary of $\Pi_h$, by $2\pi/h$-periodicity, they are the two halves of the same pick. The corresponding solutions have amplitude approximately equal to one in the physical space at $t=0$.

Denote by $u_{\pi}(x,t)$ and $u_{\eta_0}(x,t)$, $\eta_0\in(-\pi,\pi)$, the corresponding solutions of the CSE with initial data (\ref{initialdatawithoutbigrid}). Since the dispersion relation $q(\eta)=\eta^2$ for the CSE (\ref{ContSchrEqn}) is quadratic, the results above apply and the solution $u_{\pi}(x,t)$ splits into two blocks, $u_{-\pi,\gamma}^{-\pi,0}(x,t)$ and $u_{\pi,\gamma}^{0,\pi}(x,t)$, propagating along $x(t)=x^*\pm2\pi t/h$, whereas $u_{\eta_0}(x,t)$ propagates along $x(t)=x^*-2t\eta_0/h$.
\begin{figure}
  \begin{center}\includegraphics[width=5.5cm,height=3.5cm]{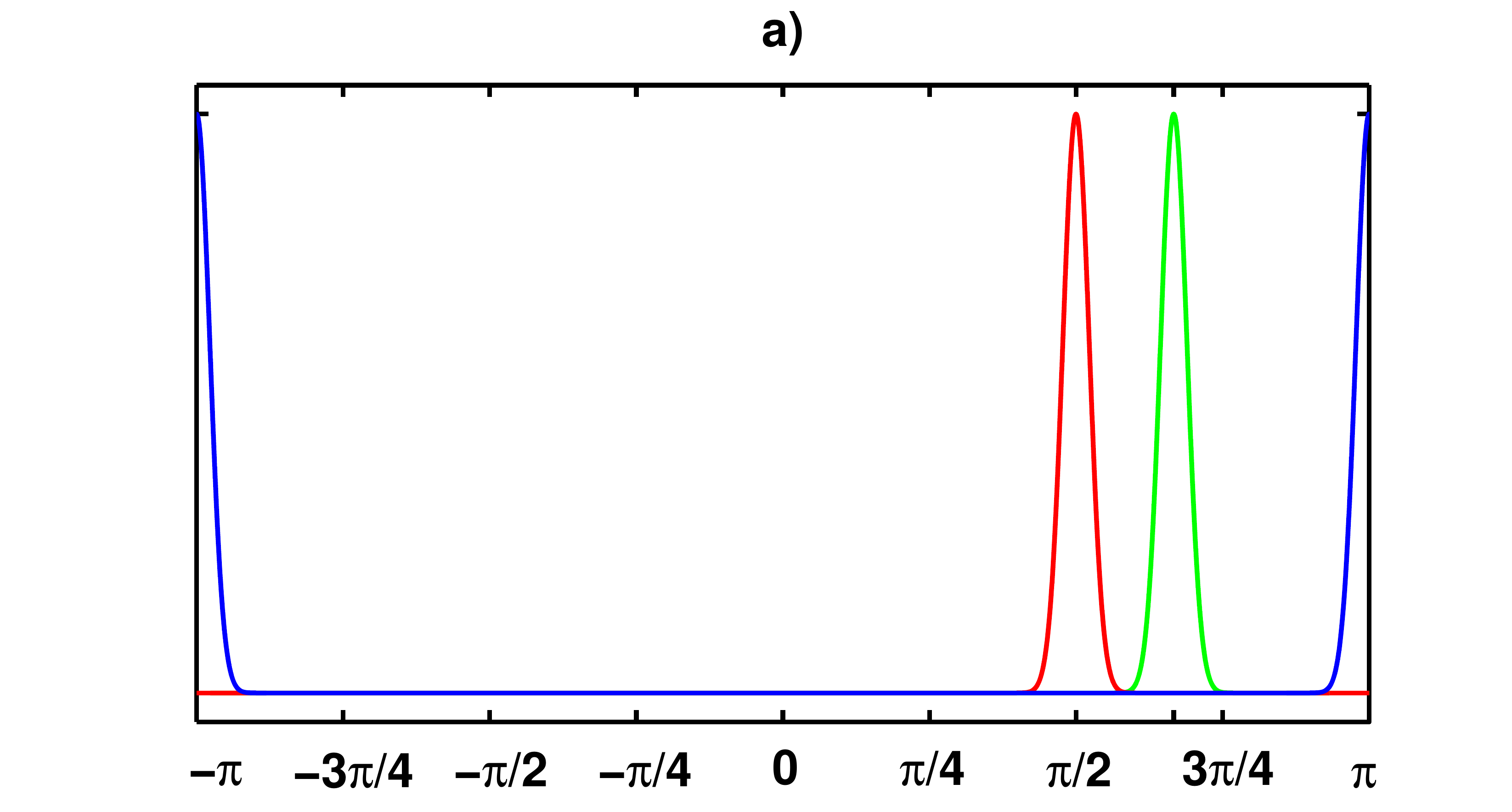}\includegraphics[width=5.5cm,height=3.5cm]{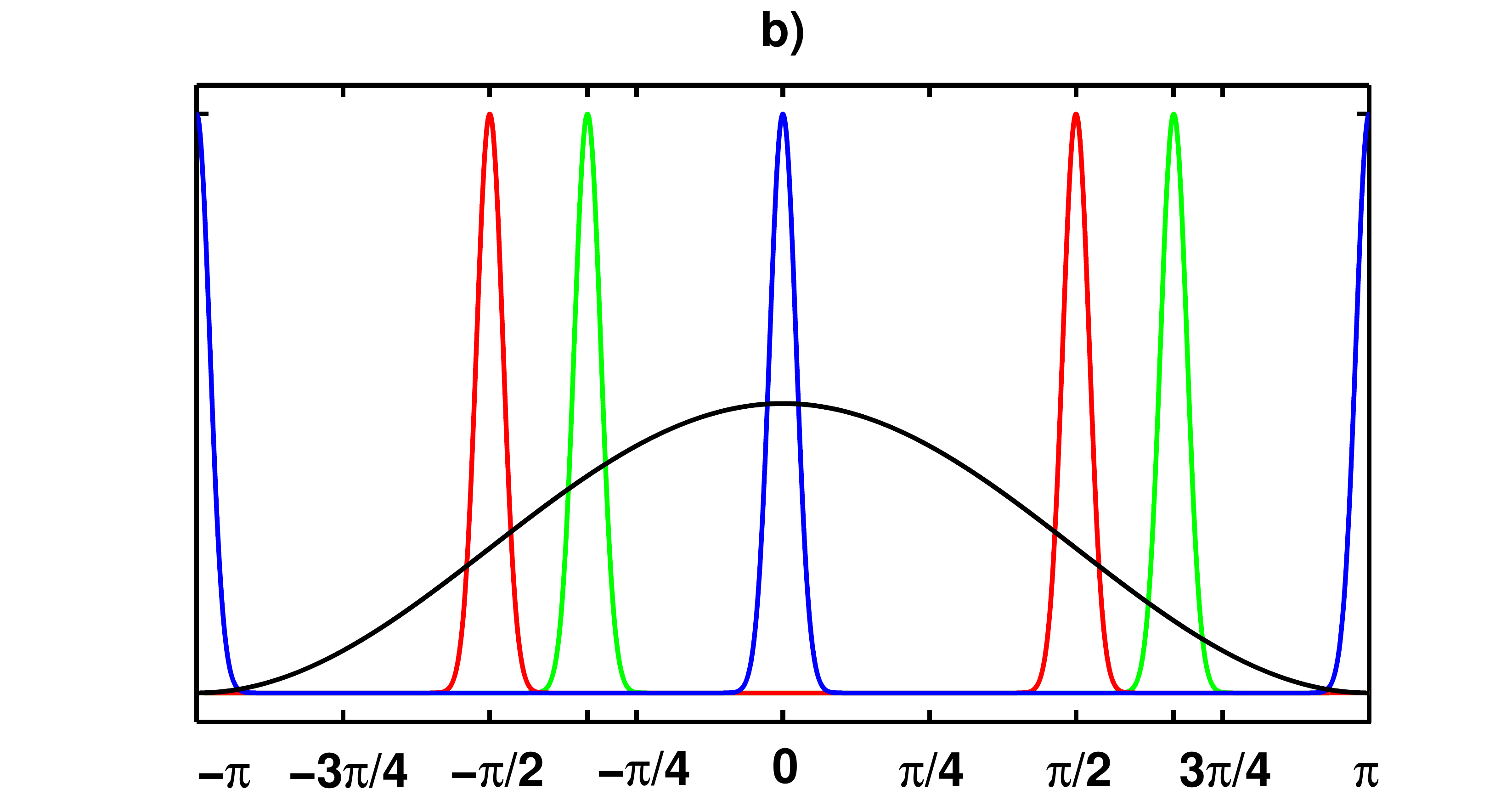}\includegraphics[width=5.5cm,height=3.5cm]{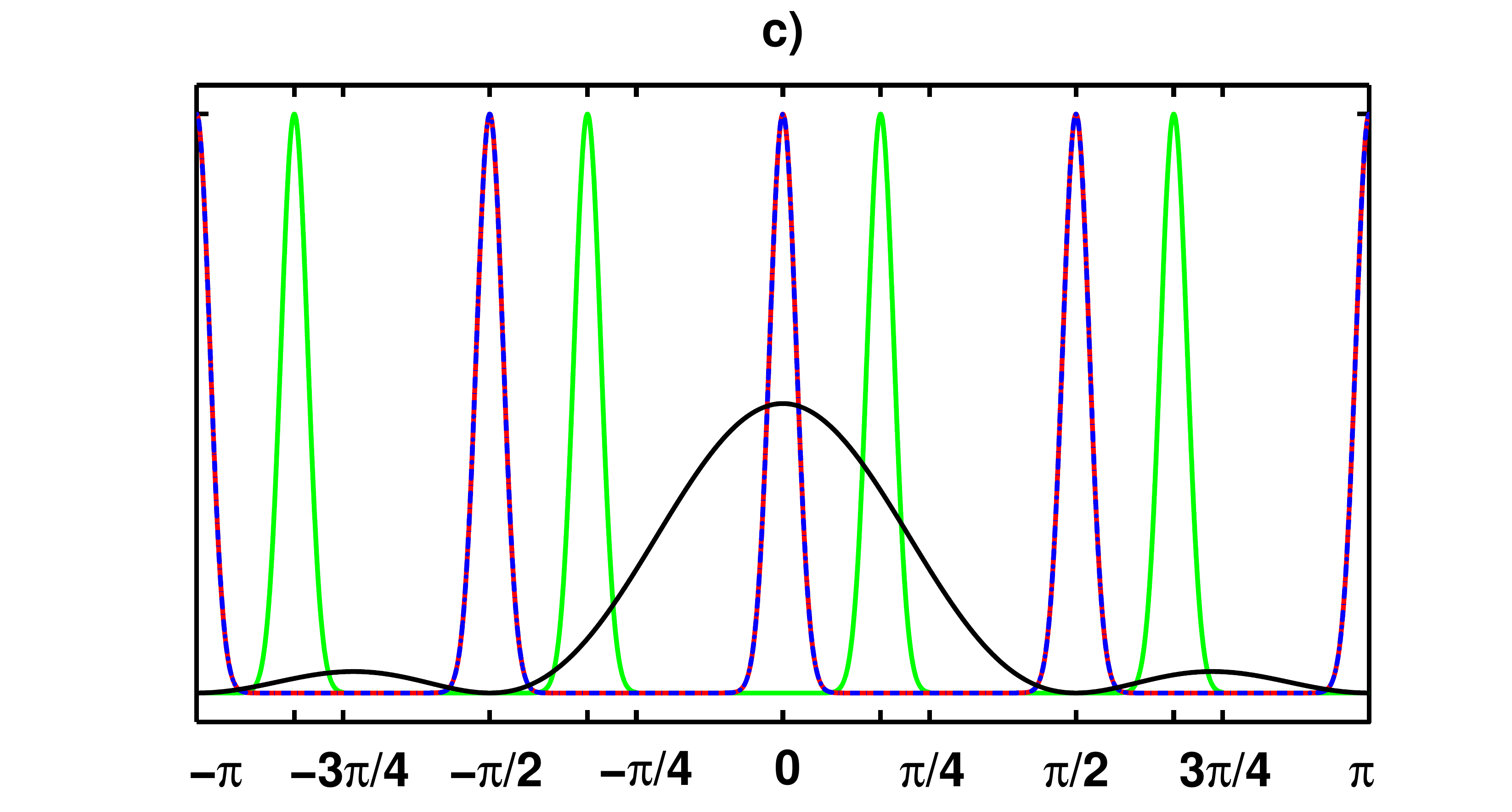}\\
  \caption{a) Initial data $\widehat{\varphi}_{\eta_0}$ with $\eta_0=\pi,\ \pi/2,\ 2\pi/3$ (blue, red, green) and their projections  $\widehat{\Lambda_k^r}^{2^kh}$ with b) $k=1$ and c) $k=2$. In black, the corresponding weights $b_k$. In c), the blue and red curves coincide, so that $\Lambda_2^{r}\varphi_{\pi}=\Lambda_2^{r}\varphi_{\pi/2}$, modulo an exponentially small error. This also means that $\varphi_{\pi}(x_{4j})$ and $\varphi_{\pi/2}(x_{4j})$ almost coincide.}\label{figr}\end{center}
\end{figure}

\begin{figure}
  \begin{center}\includegraphics[width=5.5cm,height=3.5cm]{initialdataSchrr1}\includegraphics[width=5.5cm,height=3.5cm]{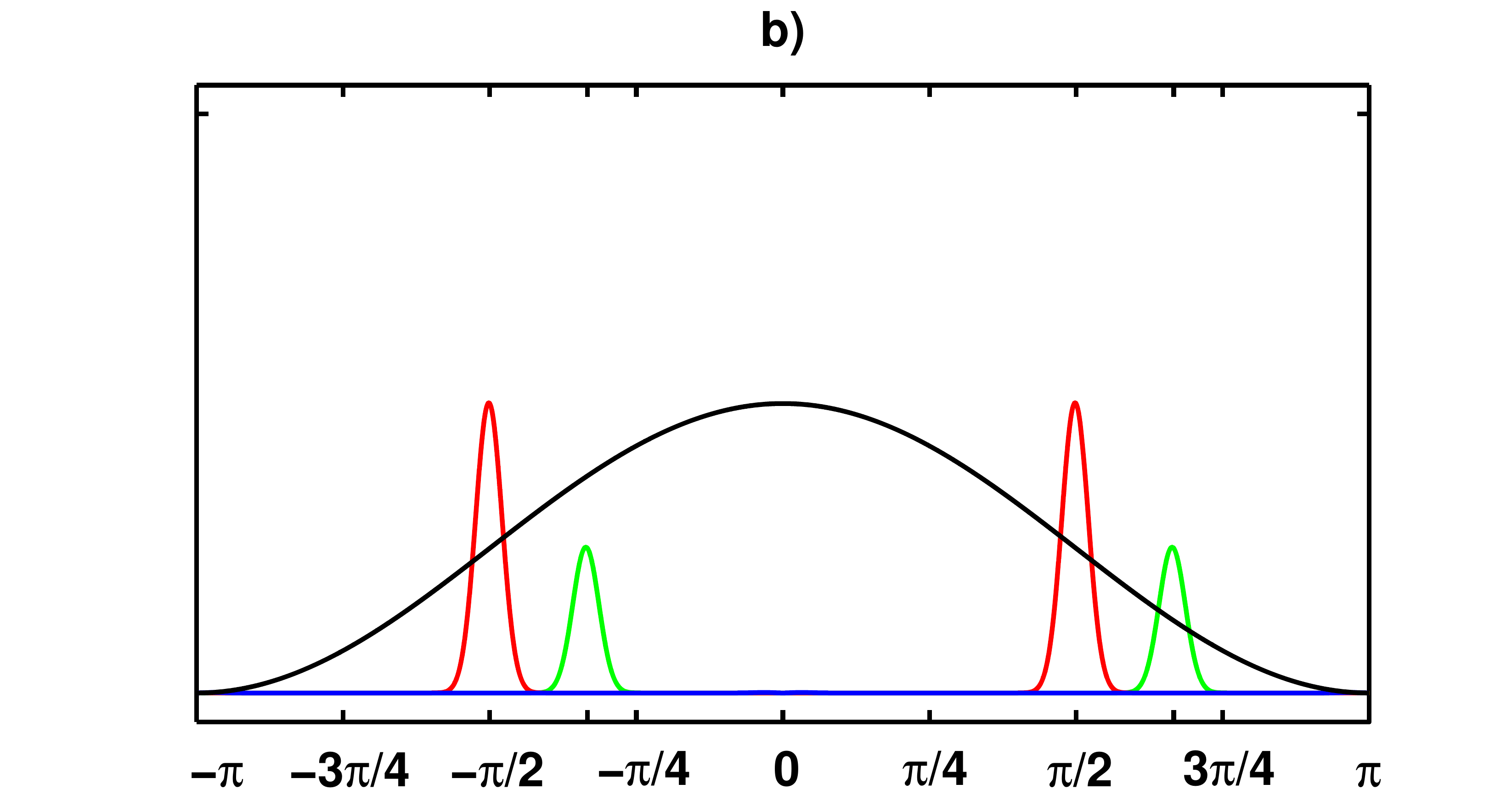}\includegraphics[width=5.5cm,height=3.5cm]{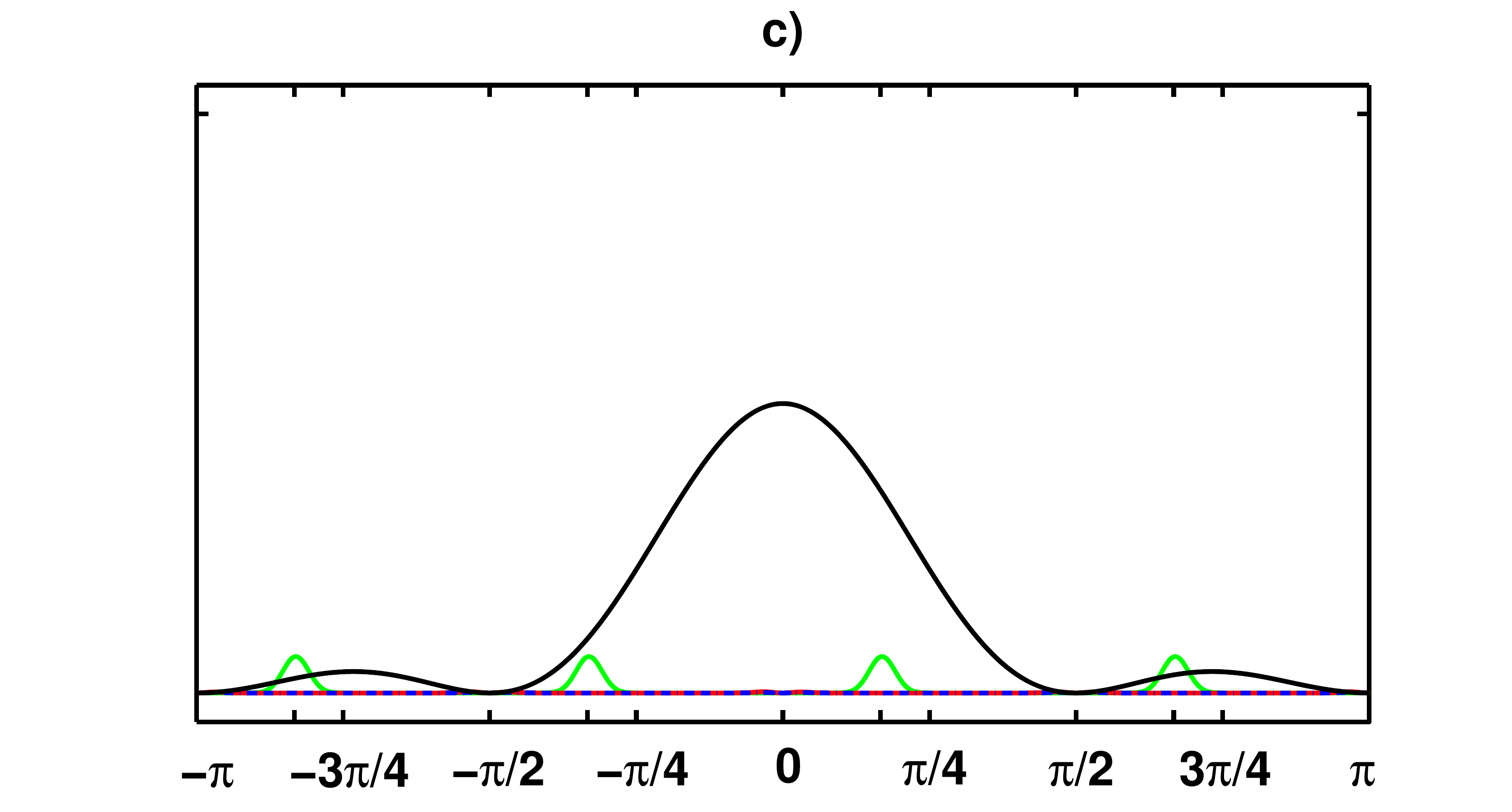}\\
  \caption{a) Initial data $\widehat{\varphi}_{\eta_0}$ with $\eta_0=\pi,\ \pi/2,\ 2\pi/3$ (blue, red, green) and their projections $\widehat{\Lambda_k^{a}}^{2^kh}$ with b) $k=1$ and c) $k=2$. In black, the corresponding weights $b_k$. The amplitude of  $\Lambda_k^{a}\varphi_{\eta_0}$ is smaller than the one of $\Lambda_k^{r}\varphi_{\eta_0}$ both in the physical and in the Fourier spaces, since $\Lambda_k^{a}\varphi_{\eta_0}\thicksim b_k(\eta_0)\Lambda_k^{r}\varphi_{\eta_0}$ and $b_k(\eta_0)<1$. Since $b_1(\pi)=0$, the blue curve in b) is almost zero. The red and green ones have amplitude $b_1(\pi/2)\sqrt{2\pi}=\sqrt{2\pi}/2$ and $b_1(2\pi/3)\sqrt{2\pi}=\sqrt{2\pi}/4$. Since $b_2(\pi)=b_2(\pi/2)=0$, the blue and red curves in c) coincide and they are almost zero. The corresponding green curve has amplitude $b_2(2\pi/3)\sqrt{2\pi}=\sqrt{2\pi}/16$.}\label{figa}\end{center}
\end{figure}

\begin{figure}
  \begin{center}\includegraphics[width=6cm,height=3.5cm]{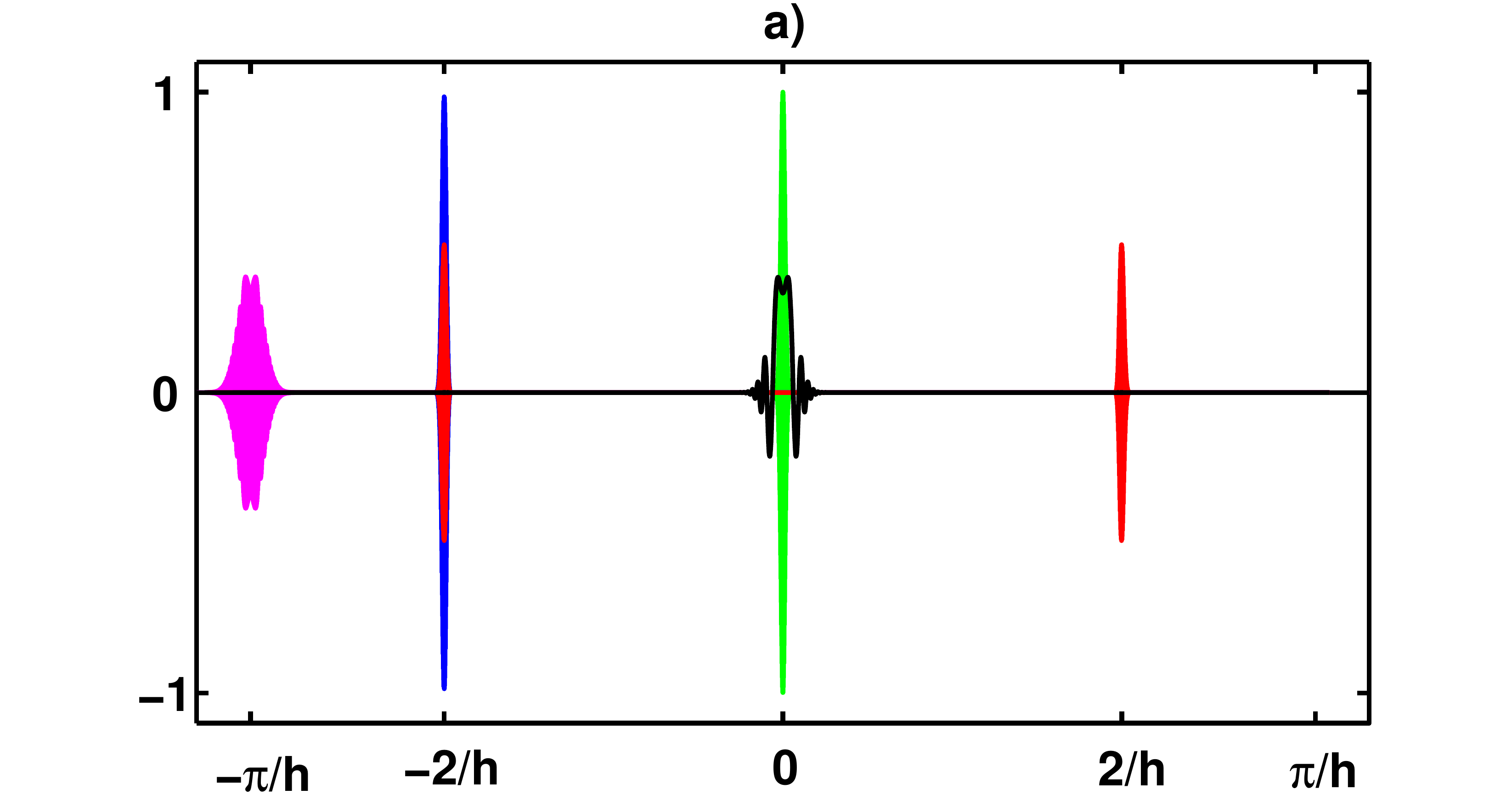}\includegraphics[width=6cm,height=3.5cm]{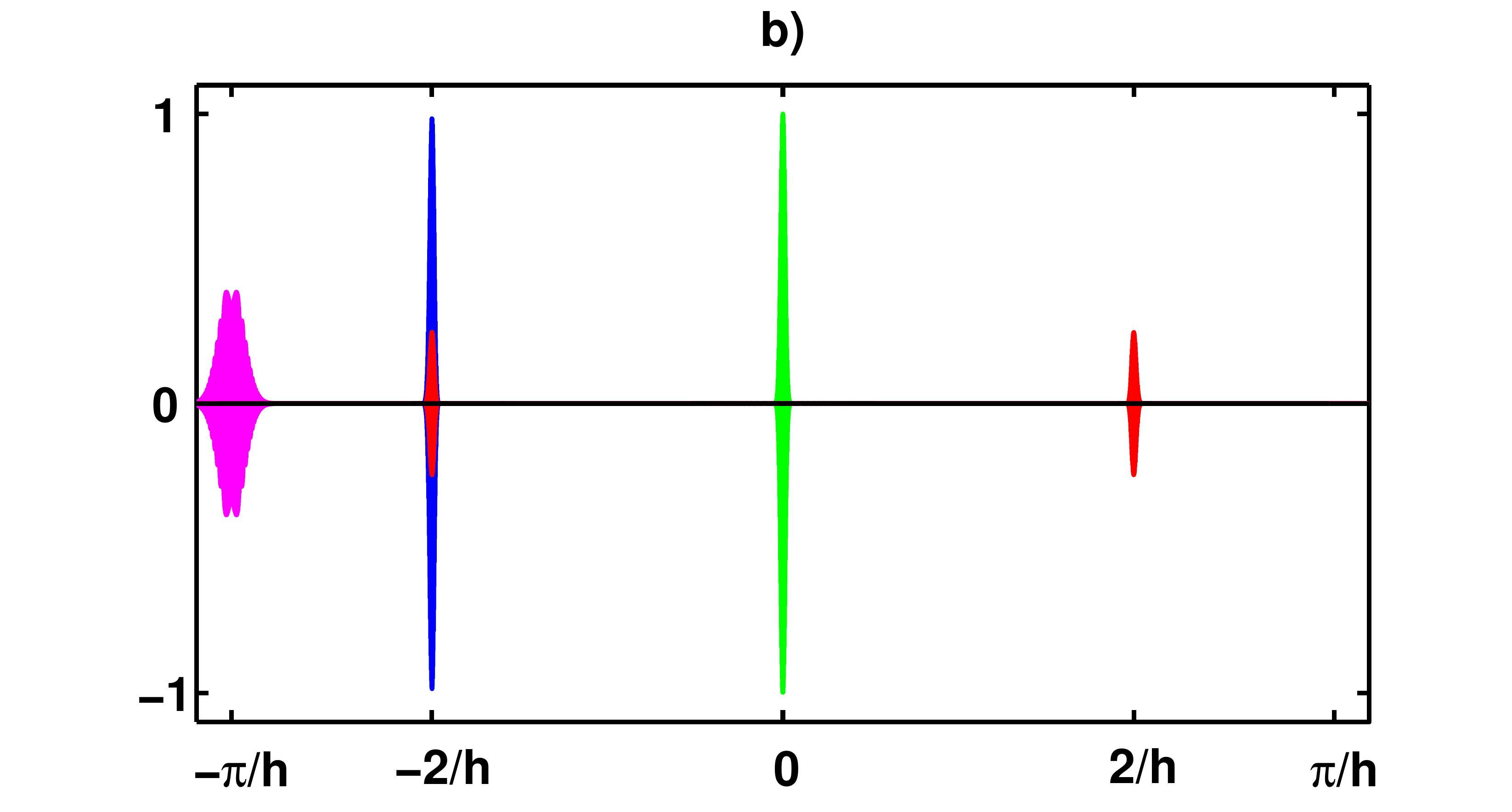}\\
\includegraphics[width=6cm,height=3.5cm]{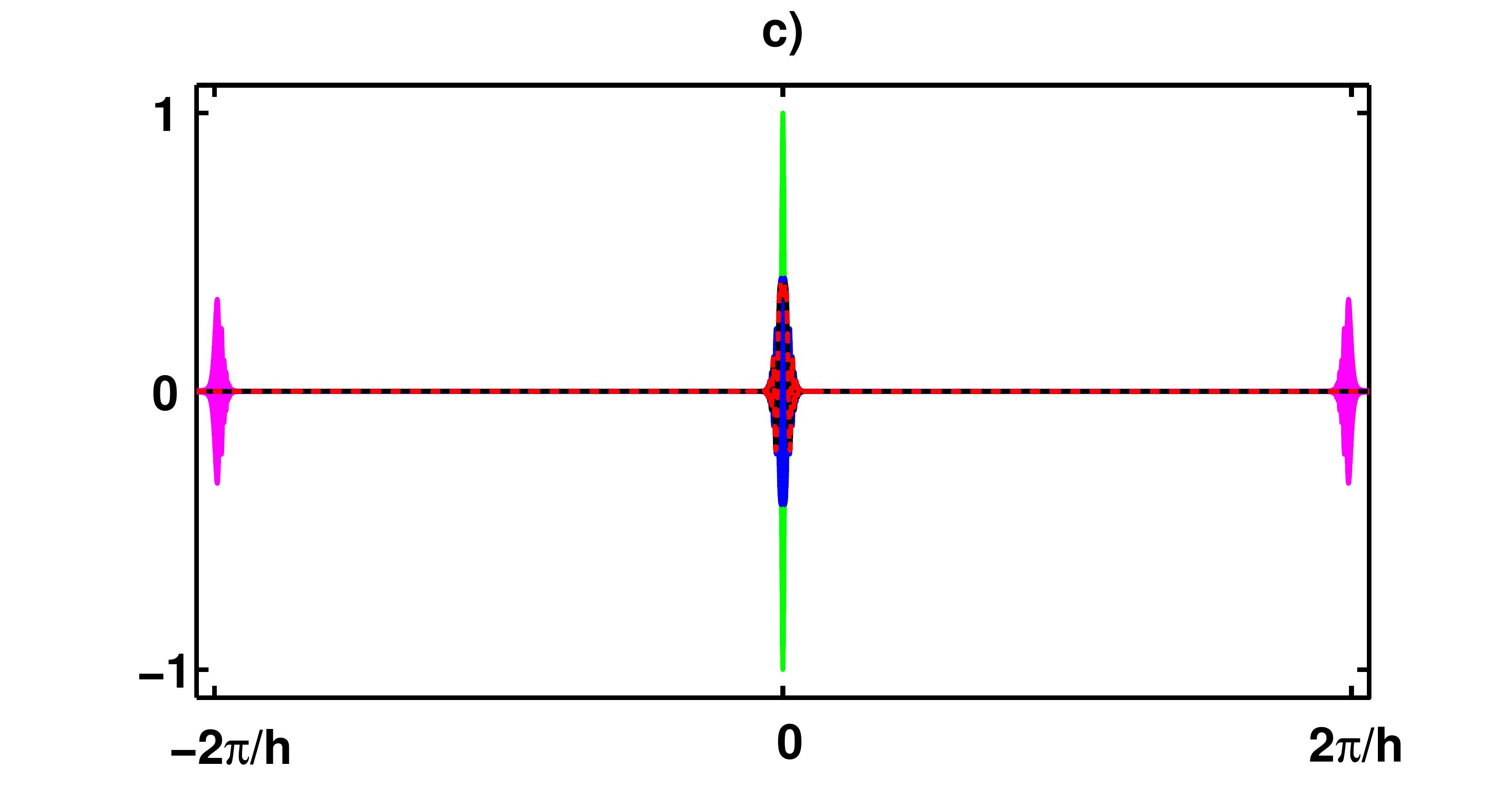}\includegraphics[width=6cm,height=3.5cm]{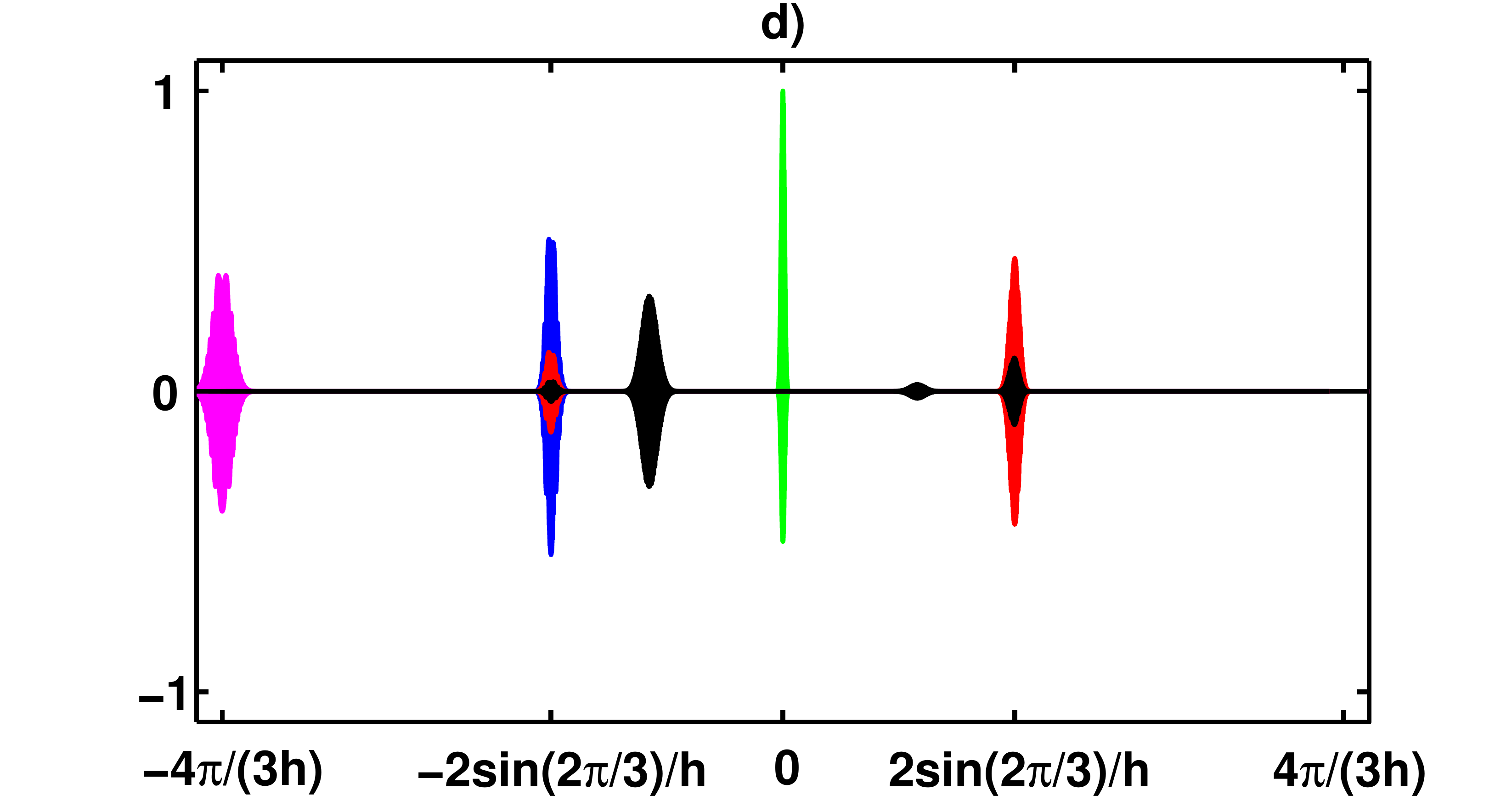}\\
  \caption{Solutions of both CSE and DSE in the physical space corresponding to the initial data $\widehat{\varphi}_{\eta_0}$ using the projections $\Lambda_k^{\alpha}$, $\alpha=r,a$, with a) $(\eta_0,\alpha)=(\pi/2,r)$, b) $(\eta_0,\alpha)=(\pi/2,a)$, c) $(\eta_0,\alpha)=(\pi,r)$ and d) $(\eta_0,\alpha)=(2\pi/3,r)$. Legend: green - solution of CSE at $t=0$, magenta - solution of CSE at $t=1$, blue - solution of DSE without filtering at $t=1$, red - solution of DSE with bi-grid of ratio $1/2$ ($k=1$) at $t=1$ and black - solution of DSE with bi-grid of ratio $1/4$ ($k=2$) at $t=1$. Description: a) the red wave packets have amplitude $b_1(\pi/2)=1/2$; b) the red wave packets have amplitude $b_1^2(\pi/2)=1/4$ and the black one has amplitude $b_2(\pi/2)=0$; c) the red and black wave packets coincide; d) in green, the initial datum is asymmetric with respect to the horizontal axis. This is an effect of the projection of the continuous initial datum into the grid. The real part of the continuous initial datum is essentially $\cos(2\pi x/(3h))\exp(-\gamma x^2/2)$, with $\gamma=h^{-1/4}$, which is symmetric, indeed. But, when sampled on $x=(jh)_{j\in\zz}$, then, for $j=0$, its amplitude is 1, but, for $j=1$, it is $\cos(2\pi/3)=-\sin(\pi/6)=-1/2$. This explains the asymmetry of the discrete plot of the datum; the red wave packets have amplitude $b_1(2\pi/3)=1/4$ and $b_1(-\pi/3)=3/4$. The black ones have amplitude $b_2(2\pi/3)=1/16$, $b_2(\pi/6)=(6+3\sqrt{3})/16$, $b_2(-5\pi/6)=(6-3\sqrt{3})/16$ and $b_2(-\pi/3)=3/16$.}\label{figpi2}\end{center}
\end{figure}

 \smallskip\noindent\textbf{Discrete solutions without filtering.} We consider the DSE with initial data as in (\ref{initialdatawithoutbigrid}). Since the dispersion relation $p(\eta)=4\sin^2(\eta/2)$ is not quadratic, we cannot apply directly the above results on quadratic dispersions relations, but this can actually be done up to a small reminder term using the Taylor expansion, as mentioned before. More precisely, for all $\eta_0\in\Pi_1$, we may split the dispersion relation as $p(\eta)=q_{\eta_0}(\eta)+r_{\eta_0}(\eta)$, where $q_{\eta_0}$ is the second-order Taylor polynomial about $\eta_0$ and $r_{\eta_0}$ is the corresponding reminder. Consider the DSE (\ref{DiscSchrEqn}) with initial data (\ref{initialdataeta0}), denote its solution by $w_{\eta_0,\gamma}^{\eta_1,\eta_2}$ and split it as $w_{\eta_0,\gamma}^{\eta_1,\eta_2}=
\widetilde{u}_{\eta_0,\gamma}^{\eta_1,\eta_2}+v_{\eta_0,\gamma}^{\eta_1,\eta_2}$, where $\widetilde{u}_{\eta_0,\gamma}^{\eta_1,\eta_2}$ is the solution corresponding to the quadratic dispersion relation $q_{\eta_0}$. It is easy to check that $\|v_{\eta_0,\gamma}^{\eta_1,\eta_2}\|_{L^2_x(\rr)}^2/\|w_{\eta_0,\gamma}^{\eta_1,\eta_2}\|_{L^2_x(\rr)}=O(th\gamma^{3/2})$, which is small in finite time intervals iff (\ref{scalegamma}) holds, so that $w_{\eta_0,\gamma}^{\eta_1,\eta_2}\thicksim\widetilde{u}_{\eta_0,\gamma}^{\eta_1,\eta_2}$.

Denote by $w_{\pi}(x,t)$ and $w_{\eta_0}(x,t)$ the solutions of the DSE with initial data (\ref{initialdatawithoutbigrid}). Since $p'(\pi)=0$, the solution $w_{\pi}$ does not propagate, illustrating the lack of uniform local smoothing effect in the discrete case (see Figure \ref{figpi2} - c), the blue curve). For $\eta_0=\pi/2$, since $p''(\pi/2)=0$, the spatial support of $w_{\pi/2}$ does not expand as time evolves, which agrees with the lack of uniform $L^q_t\ell^p_x$-integrability properties (see Figure \ref{figpi2} - a) and b), the blue curve).

\smallskip\noindent\textbf{Bi-grid solutions.} For the DSE, we now consider initial data obtained by firstly projecting the data   (\ref{initialdatawithoutbigrid}) from the fine grid to the coarse one by one of the two projections (\ref{projectionsbigrid}) and then extending those projections by linear interpolation from the grid of size $2^kh$ to the one of size $h$. We identify several different cases:

\smallskip\noindent\textit{A) $\eta_0=\pi$, $\alpha=r$.} In the Fourier space, the projection $\Lambda_k^r\varphi_{\pi}$ has picks at $2l\pi/2^kh$, for $0\leq |l|\leq 2^{k-1}$. The picks at $2l\pi/2^kh$, $1\leq |l|\leq 2^{k-1}$, are canceled as $h\to 0$ since the weight $b_k$ vanishes exactly at those points. The wave packet for $l=0$, being located at $\xi=0$, \textit{does not propagate}, but it \textit{does not oscillate} either and the discrete smoothing property holds uniformly as $h\to 0$ (see Figure \ref{figpi2} - c), the black curve, for $(\eta_0,k)=(\pi,2)$). This packet also decays with the rate of the CSE since both dispersion relations (continuous and discrete one) are tangent at $\xi=0$.

\smallskip\noindent\textit{B) $\eta_0\in((2l^*-1)\pi/2^k,(2l^*+1)\pi/2^k)$, for some $-2^{k-1}+1\leq l^*\leq2^{k-1}-1$, $\alpha=r$.} Set  $\eta_0^*:=\eta_0-2l^*\pi/2^k\in\Pi_{2^kh}$. Similarly, if $\eta_0\in(-\pi,-(2^k-1)\pi/2^k)\cup((2^k-1)\pi/2^k,\pi)$, set $\eta_0^*:=\eta_0\pm\pi$. There are two cases: i) $\eta_0^*=0$. Then the projection $\Lambda_k^r\varphi_{\eta_0}$ has picks in the Fourier space at $\xi=2l\pi/(2^kh)$, with $0\leq|l|\leq2^{k-1}$ and the analysis follows the one in part A) (see Figure \ref{figpi2} - a), for $(\eta_0=\pi/2,k=2)$); ii) $\eta_0^*\in(-\pi/2^k,\pi/2^k)\setminus\{0\}$. Then the projection $\Lambda_k^r\varphi_{\eta_0}$ has picks in the Fourier space at $(\eta_0^*+2l\pi/2^k)/h$, where $-2^{k-1}+1-s\leq l\leq 2^{k-1}-s$ and $s=(1+\mbox{sign}(\eta_0^*))/2$. The solution $w_{\eta_0,k}^{r}$ is a superposition of $2^k$ wave packets propagating along the lines $x(t)=x^*-tv_l$  with velocity $v_l=2\sin(\eta_0^*+2l\pi/2^k)/h$. Observe that $v_l=-v_{-2^{k-1}+l}$, $1-s\leq l\leq 2^{k-1}-s$, i.e. there are $2^{k-1}$ pairs of wave packets constituted by one wave packet going in each direction with the same velocity, and $v_{1-s}<v_{2^{k-1}-s}<v_{2-s}<v_{2^{k-1}-1-s}<\cdots<v_{2^{k-2}-s}<v_{2^{k-2}+1-s}$. If $\mbox{sign}(\eta_0)\not=\mbox{sign}(\eta_0^*)$, then the wave packet of largest amplitude changes the direction with respect to the solution without filtering (see Figure \ref{figpi2} - d), the red curve, for $(\eta_0,k)=(2\pi/3,1)$). In this case, the dispersive properties are verified uniformly.

\smallskip\noindent \textit{C) $\eta_0=(2l^*+1)/2^k$, for some $-2^{k-1}\leq l^*\leq 2^{k-1}-1$, $\alpha=r$.} Then
the projection $\Lambda_k^r\varphi_{\eta_0}$ has picks in the Fourier space at $\xi=(2l+1)\pi/(2^kh)$, for all $-2^{k-1}\leq l\leq 2^{k-1}-1$. The solution $w_{\eta_0,k}^{r}$ is a superposition of $2^k$ wave packets, each one propagating along the line $x(t)=x^*-tv_l$ with velocity $v_l=2\sin((2l-1)\pi/2^k)/h$. With respect to B) - ii., we have $v_l=v_{2^{k-1}+1-l}$ for all $1\leq l\leq 2^{k-2}$, i.e. the $2^{k-1}$ wave packets going to the left can be grouped into $2^{k-2}$ pairs constituted by two wave packets having the same velocity which collapse at any time $t>0$, and $v_1<v_2<\cdots<v_{2^{k-2}}$. For $(\eta_0,k)=(\pi/2,1)$, the two blocks in the solution do not spread as time evolves (see Figure \ref{figpi2} - a), the red curve), which confirms the fact that the bi-grid algorithm with ratio $1/2$ is not enough to reestablish the uniformity of the dispersive estimates, as predicted by the theory in \cite{LivEZDispSchr}.

\smallskip\noindent \textit{D) $\alpha=a$, $\eta_0\in(-\pi,\pi]$.} Remark that $\Lambda_k^{a}\varphi_{\eta_0}\thicksim b_k(\eta_0)\Lambda_k^{r}\varphi_{\eta_0}$. Then: i) the solutions corresponding to $\Lambda_k^{a}$ are of smaller amplitude than those corresponding to $\Lambda_k^{r}$, because they involve the factor $b_k(\eta_0)<1$ (see the red curves in Figure \ref{figpi2} - a) and b)); ii) when $\eta_0$ is such that $b_k(\eta_0)=0$, the projection $\Lambda_k^{a}$ is too strong, in the sense that the corresponding numerical solutions tends to zero as $h\to0$ (see Figure \ref{figa} - b), the blue curve for $(\eta_0,k)=(\pi,1)$, or c), the blue/red curves, for $(\eta_0,k)=(\pi,2),(\pi/2,2)$. Also the black curve in Figure \ref{figpi2} - a), for $(\eta_0,k)=(\pi/2,2)$). The cancellations can be explained by the fact that for plane waves of the form $\varphi_j=\exp(i\eta_0 j)$, $j\in\zz$, the identity $\Lambda_k^{a}\overrightarrow{\varphi}^h\equiv0$ holds.

\smallskip\noindent\textbf{Conclusions.} In this article, some subtle phenomena related to the dispersivity (and the lack of) have been described in the context of the Schr\"{o}dinger equation and its numerical approximation schemes. We have shown that, choosing appropriate high frequency wave packets, one may confirm the predictions of the theory (cf. \cite{LivEZDispSchr}) and gain new insight about the complex behavior of such solutions which, depending on the high-frequency around which they concentrate, may or not propagate, or propagate but not disperse. We also analyze and describe how these solutions are affected by the application of bi-grid filtering techniques. The bi-grid solutions can be decomposed into several wave packets moving at different velocities. Some of them are attenuated by the bi-grid weights or by the way one projects the initial data from the fine grid to the coarse one. In this way, the bi-grid mechanisms may yield non-oscillatory, smooth solutions or solutions vanishing as $h\to0$, even if the original initial data do not have those behaviors. Our numerical simulations confirm the predictions of the theory in \cite{LivEZDispSchr}, in the sense that the bi-gird technique, implemented with mesh-ratio $1/4$, ensures the dispersive properties of numerical solutions.

\smallskip\noindent\textbf{Open problems.} There are several issues, related to those treated in this paper, worth to be explored: i) So far we have plotted solutions based on their Fourier description. It would be interesting to do the same, but using numerical approximations schemes. But then one has to take into account the changes that the time discretizations (cf. \cite{LivEZDispSchr}) or some kind of boundary conditions, for example the transparent ones (cf. \cite{TransparentBCSchr}), could introduce in the dispersive properties; ii) the relevance of these high frequency phenomena is still to be explored in the context of the non-linear Schr\"{o}dinger equations. The theoretical results in \cite{LivEZConv} show that dispersive schemes ensure a polynomial convergence order which improves the logarithmic one one gets for standard finite-difference schemes by energy methods. Also the effect of splitting methods (\cite{LivSplit}, \cite{FaouSplit}) on these high frequency wave packets is worth investigating; iii) The issues addressed in this Note are totally open for non-uniform grids or Schr\"{o}dinger equations in heterogeneous media.

\smallskip\noindent\textbf{Acknowledgements.} Both authors were partially supported by the Grant MTM2008-03541 of the MICINN, Spain, project PI2010-04 of the Basque Government and the ERC Advanced Grant FP7-246775 NUMERIWAVES.

\end{document}